\documentclass[12pt,a4paper]{article}
\usepackage{amsmath,amsthm}
\usepackage{amsfonts}
\newtheorem{theorem}{Theorem}[section]
\newtheorem{lemma}[theorem]{Lemma}
\newtheorem{corollary}[theorem]{Corollary} 
\newtheorem{proposition}[theorem]{Proposition}

\title{\bf On  subgroups in division rings of type 2}

\author{Bui Xuan Hai\footnote{Corresponding author, Faculty of Mathematics and Computer Science, University of Science, VNU-HCM, 227 Nguyen Van Cu Str., Dist. 5, HCM-City, Vietnam,  e-mail: bxhai@hcmus.edu.vn}, Trinh Thanh Deo\footnote{Faculty of Mathematics and Computer Science, University of Science, VNU-HCM, 227 Nguyen Van Cu Str., Dist. 5, HCM-City, Vietnam,  e-mail: ttdeo@hcmus.edu.vn}, and Mai Hoang Bien\footnote{Department of Basic Sciences, University of Architecture, 196 Pasteur Str., Dist. 1, HCM-City, Vietnam, e-mail:  maihoangbien012@yahoo.com}}
\date{}
\begin{document}
\baselineskip=18pt
\maketitle

\newcommand{\dpcm}{ \hfill \rule{3mm}{3mm}}
\def\Q{\mathbb{Q}}
\def\F{\mathbb{F}}
\newcommand{\ts}[1]{\langle #1\rangle}

\begin{abstract} 
\baselineskip=18pt
Let $D$ be a division ring with center $F$. We say that $D$ is a {\em division
ring of type $2$} if for every two elements $x, y\in D,$ the division subring
$F(x, y)$ is a finite dimensional vector space over $F$. In this paper we
investigate multiplicative subgroups in such a ring.
\end{abstract}



\newpage

\section{Introduction}

In the theory of division rings, one of the problems is to determine which
groups can occur as multiplicative groups of non-commutative division rings. 
There are some interesting results relating to this problem. Among them  we
note the famous discovery of  Wedderburn in 1905, which states that  {\it if
$D^*$ is a finite group, then $D$ is commutative}, where $D^*$ denotes the
multiplicative group of $D$. Later, L. K. Hua (see, for example, in [12,
p. 223]) proved that the multiplicative group of a non-commutative division
ring cannot be solvable. Recently, in \cite{hai-thin2} it was shown that the
group $D^*$ cannot even be locally nilpotent. Note also Kaplansky's Theorem
(see [12,(15.15), p. 259]) which states that if the group $D^*/F^*$ is
torsion, then $D$ is commutative, where  $F$ is the center of $D$. 
Some other results of this kind can be found for example, in
\cite{Ak1}-\cite{Ak3}, \cite{hai-huynh}-\cite{hai-thin2},...

In this paper we consider this question for division rings of type $2$.
Recall that a division ring $D$ with center $F$ is said to be {\em division
ring of type $2$} if for every two elements $x, y\in D,$ the division 
subring $F(x, y)$ is a finite dimensional vector space over $F$. This
concept is an extension of that of locally finite division rings. By
definition, a division ring $D$ is {\em centrally finite} if it is a finite
dimensional vector space over its center $F$ and $D$ is {\em locally finite}
if for every finite subset $S$ of $D$, the division subring $F(S)$ 
generated by $S\cup F$ in $D$ is a finite dimensional vector space over
$F$. There exist locally finite division rings which are not centrally finite
(it is not hard to give some examples). Of course, every locally finite
division ring is a ring of type $2$. However, at present no example of a
division ring of type $2$ is known which is not locally finite. The
difficulties are related with the following famous longstanding conjecture
known as the Kurosh Problem for division rings \cite{kha}. Recall that a
division ring $D$ is {\em algebraic} over its center $F$ (briefly, $D$ is 
{\em algebraic}), if every element of $D$ is algebraic over $F$. 
Clearly, a locally finite division ring is algebraic. Kurosh conjectured that
any algebraic division ring is locally  finite. Unfortunately, this problem
remains still unsolved in general, it is answered in the affirmative for the
following special cases: for $F$ uncountable \cite{row}, $F$ finite
\cite{lam}, and for $F$ having only finite algebraic field extensions (in
particular for $F$ algebraically closed). The last case follows from the
Levitzki-Shirshov Theorem which states that {\em any algebraic algebra of
bounded degree is locally finite} (see e.g. \cite{dren}, \cite{kha}). The
answer for the case of finite $F$ is due to Jacobson who proved that {\em an
algebraic division ring $D$ is commutative provided its center is finite}
(see, for example, \cite{lam}). Later, more general theorems of this kind
(known as commutativity theorems) were proved by Jacobson and Herstein. For
more information we refer to [9, Ch. 3]. Finally, we would like to note that
the results obtained in this paper for division rings of type $2$ have not
been proved elsewhere before for locally finite division rings. So, at least
(in the fortunate case if the Kurosh Problem will be answered in the
affirmative, as we would like to see) our results generalize previous results
for the finite dimensional case. 

Throughout this paper the following notation will be used consistently: $D$
denotes a division ring with center $F$ and $D^*$ is the multiplicative group
of $D$. If $S$ is a nonempty subset of $D$, then we  denote  by $F[S]$ and
$F(S)$ the subring and the division subring of $D$ generated by $S$ over $F$,
respectively. The symbol $D'$ is used to denote the derived group 
$[D^*, D^*]$. An element $x$ in $D$ is said to be {\em radical} over a subring
$K$ of $D$ if there exists some positive integer $n(x)$ depending on $x$ such
that $x^{n(x)}\in K$. A nonempty subset $S$ of $D$ is {\em radical} over $K$ if
every element of $S$ is radical over $K$. We denote by $N_{D/F}$ and
$RN_{D/F}$ the norm and the reduced norm, respectively. Finally, if $G$ is any
group then we always use the symbol $Z(G)$ to denote the center of $G$.

\section{Finitely generated subgroups}

The main purpose in this section is to prove that in any non-commutative
division ring  of type $2$  there are no finitely generated subgroups
containing the center. 

\begin{lemma}\label{lem:2.1}
Let $D$ be a division ring with center $F$, $D_1$ be a division subring of
$D$ containing $F$. Suppose that $D_1$ is a finite dimensional vector space
over $F$ and $a\in D_1$. Then, $N_{D_1/F}(a)$ is periodic if and only if
$N_{F(a)/F}(a)$ is periodic.  
\end{lemma}

\begin{proof}
Let $F_1=Z(D_1)\supset F$, $m^2=[D_1:F_1]$ and $n=[F_1(a):F_1]$.
By [4, Lemma 3, p.145] and [4, Corollary 4, p. 150], we have
$$N_{D_1/F_1}(a)=[RN_{D_1/F_1}(a)]^m=[N_{F_1(a)/F_1}(a)]^{m^2/n}.$$
Now, using the Tower formulae for the norm (cf. \cite{dra}), from the
equality above we get
$$N_{D_1/F}(a)=[N_{F_1(a)/F}(a)]^{m^2/n}.$$
Since $a\in F(a), we have N_{F_1(a)/F(a)}(a)=a^k$, where
$k=[F_1(a):F(a)]$. Therefore 
$$N_{F(a)/F}(a^{k})=N_{F(a)/F}(N_{F_1(a)/F(a)}(a))=N_{F_1(a)/F}(a).$$
It follows that $N_{D_1/F}(a)=[N_{F(a)/F}(a)]^{km^2/n}$,
and the conclusion is now obvious.
\end{proof}

The following proposition is useful. In particular, it is needed to prove the
subsequent theorem. 

\begin{proposition}\label{prop:2.2}
Let $D$ be a division ring with center $F$. If $N$ is a subnormal subgroup of
$D^*$ then $Z(N)=N\cap F$. 
\end{proposition}

\begin{proof}
If $N$ is contained in $F$ then there is nothing to prove. Thus, suppose that
$N$ is non-central. By [15, 14.4.2, p. 439], $C_D(N)=F$. Hence 
$Z(N)\subseteq N\cap F$. Since the inclusion $N\cap F\subseteq Z(N)$ is
obvious, $Z(N)= N\cap F$. 
\end{proof}

\begin{theorem}\label{thm:2.3}
Let $D$ be a division ring of type $2$. Then  $Z(D')$ is a torsion group.

\end{theorem}

\begin{proof}
By Proposition \ref{prop:2.2}, $Z(D')=D'\cap F^*$. Any element $a\in Z(D')$
can be written in the form $a=c_1 c_2\ldots c_r$, where $c_i=[x_i, y_i]$ with
$x_i, y_i\in D^*$ for $i\in\{1, \ldots, r\}$. Put $D_1=D_2:= F(c_1, c_2)$,
$D_3:= F(c_1 c_2, c_3)$, $\ldots$, $D_r:= F(c_1...c_{r-1}, c_r)$ and
$F_i=Z(D_i)$ for $i\in\{1, \ldots, r\}$. Since $D$ is of type $2$,
$[D_i:F]<\infty$.

Since $N_{F(x_i, y_i)/F}(c_i)=1$, by Lemma \ref{lem:2.1}, $N_{F(c_i)/F}(c_i)$
is periodic. Again by Lemma~\ref{lem:2.1}, $N_{D_i/F}(c_i)$ is
periodic. Therefore, there exists some positive integer $n_i$ such that 
$N_{D_i/F}(c_i^{n_i})=1$. Recall that $D_2=D_1$. Hence we get
$$N_{D_2/F}(c_1 c_2)^m=N_{D_2/F} (c_1)^m  N_{D_2/F}(c_2)^m=1,$$
where $m=n_1 n_2$. Again by Lemma \ref{lem:2.1}, $N_{F(c_1 c_2)/F}(c_1 c_2)$
is periodic; hence  $N_{D_3/F}(c_1 c_2)$ is periodic. By induction,
$N_{D_r/F}(c_1... c_{r-1})$ is periodic. Suppose that
$N_{D_r/F}(c_1... c_{r-1})^n=1$. Then 
$$N_{D_r/F}(a^n)=N_{D_r/F}(c_1... c_{r-1})^n N_{D_r/F}(c_r)^n=1.$$
Hence, $a^{n[D_r:F]}=1$.
Therefore, $a$ is periodic. Thus $Z(D')$ is torsion.
\end{proof}

\begin{corollary}\label{cor:2.4}
Let $D$ be a non-commutative ring of type $2$ with center $F$. Then
$D'\setminus Z(D')$ contains no elements purely inseparable over $F$. 
\end{corollary}

\begin{proof} Suppose that $a\in D'\setminus Z(D')$ is purely inseparable over
$F$. Then, there exists some positive integer $m$ such that 
$a^{p^m}\in F$. Since $Z(D')=D'\cap F$ (by Proposition 2.2), 
$a^{p^m}\in Z(D')$. By Theorem \ref{thm:2.3}, there exists some positive
integer $r$ such that $a^{rp^m}=1$. Denote by $k$ the order of $a$ in the
group $D^*$. If $p$ divides $k$, then $k=pt$ and we have 
$$1=a^k=a^{pt}=(a^t)^p.$$

Consequently, $a^t=1$, which is impossible in view of the choice of $k$. Now,
suppose that $p$ does not divide $k$. Then, $(k, p^m)=1$ and $\alpha k+\beta
p^m=1$ for some integers $\alpha$ and $\beta$. Therefore, we have 
$$a=a^{\alpha k+\beta p^m}=(a^k)^\alpha.a^{\beta p^m}=(a^{p^m})^\beta\in F.$$

Consequently, $a\in F\cap D'=Z(D')$, a contradiction.
\end{proof}

Note that in \cite{mah} the author proved that $Z(D')$ is finite if $D$ is
centrally finite. In virtue of this fact, he expressed his ideas  that
$Z(D')$ is torsion for any division ring $D$ algebraic over its center, but he
has not been able to prove this. Therefore, Theorem~\ref{thm:2.3} represents
some progress in this direction. Moreover (and this is more important for our
purpose), we need this theorem to establish the main result in the present
section. In fact, we shall prove that in a division ring $D$ of type $2$ with
center $F$, there are no finitely generated subgroups containing
$F^*$. Consequently, if $D$ is of type $2$ and $D^*$ is finitely generated,
then $D$ is a  field. Note that if the multiplicative group of a field is
finitely generated, then it is finite. So, if $D$ is of type $2$ and $D^*$ is
finitely generated, then $D$ is even a finite field. Our next theorem strongly
generalizes the result obtained in [2, Theorem 1] which states that, if $D$ is
centrally finite and $D^*$ is finitely generated, then $D$ is commutative.
 
\begin{theorem}\label{thm:2.6}
Let $D$ be a non-commutative division ring of type $2$ with center $F$ and
suppose that $N$ is a subgroup of $D^*$ containing $F^*$. Then $N$ is not
finitely generated. 
\end{theorem}

\begin{proof}
Suppose that there is a finitely generated subgroup 
$N=\ts{x_1, x_2, \ldots, x_n}$ of $D^*$ containing $F^*$. Then, in virtue of
[15, 5.5.8, p. 113], $F^*N'/N'$ is a finitely generated abelian group, where
$N'$ denotes the derived subgroup of $N$. 

\noindent
{\em Case 1: $char(D)=0$.}

Then, $F$ contains the field $\Q$ of rational numbers and it follows that
$\Q^*/(\Q^*\cap N')\simeq \Q^*N'/N'$. Since $F^*N'/N'$ is finitely generated,
$\Q^*N'/N'$ is finitely generated and consequently   $\Q^*/(\Q^*\cap N')$ is
finitely generated. Consider an arbitrary element $a\in \Q^*\cap N'$. Then
$a\in F^*\cap D'=Z(D')$. By Theorem \ref{thm:2.3}, $a$ is periodic.
Since $a\in \Q$, we get $a=\pm{1}$. Thus, $\Q^*\cap N'$ is finite. Since
$\Q^*/(\Q^*\cap N')$ is finitely generated, $\Q^*$ is finitely generated,
which is impossible. 

\noindent
{\em Case 2: $char(D)=p > 0$.}

Denoting by $\F_p$ the prime subfield of $F$, we shall prove that $F$ is
algebraic over $\F_p$. In fact, suppose that $u\in F$ and $u$  is
transcendental over $\F_p$. Then, the group $\F_p(u)^*/(\F_p(u)^*\cap N')$
considered as a subgroup of $F^*N'/N'$ is finitely generated. Consider an
arbitrary element $f(u)/g(u)\in \F_p(u)^*\cap N'$, where $f(X), g(X)\in
\F_p[X], ((f(X), g(X))=1$ and $g(u)\neq 0$. As above, we have
$f(u)^s/g(u)^s=1$ for some positive integer $s$. Since $u$ is transcendental
over $\F_p$, it follows that $f(u)/g(u)\in \F_p$. Therefore, $\F_p(u)^*\cap
N'$ is finite and consequently, $\F_p(u)^*$ is finitely generated, so 
$\F_p(u)$ is finite field, which is impossible. Hence $F$ is algebraic over
$\F_p$ and it follows that $D$ is algebraic over $\F_p$. Now, in virtue of
Jacobson's Theorem [12, (13.11), p. 219], $D$ is commutative, a contradiction.
\end{proof}

\begin{corollary}\label{cor:2.7}
Let $D$ be a division ring of type $2$. If $D^*$ is finitely generated, then
$D$ is a finite field. 
\end{corollary}

If $M$ is a finitely generated maximal subgroup of $D^*$, then clearly $D^*$
is finitely generated. So, the next result follows immediately from
Corollary \ref{cor:2.7}. 

\begin{corollary}\label{cor:2.8}
Assume that  $D$ is  a division ring of type $2$. If  $D^*$  has a  finitely
generated maximal subgroup, then $D$ is a finite field. 
\end{corollary}

In the same way as in the proof of Theorem \ref{thm:2.6}, we obtain the
following corollary. 

\begin{corollary}\label{cor:2.9}
Assume that $D$ is a non-commutative division ring of type $2$ with center
$F$ and $S$ is a subgroup of $D^*$. If $N=SF^*$, then $N/N'$ is not finitely
generated. 
\end{corollary}

\begin{proof}
Suppose that $N/N'$ is finitely generated. Since $N'=S'$ and 
$F^*/(F^*\cap S') \simeq S'F^*/S'$, it follows that $F^*/(F^*\cap S')$ is a
finitely generated abelian group. Now, in the same way as in the proof of
Theorem \ref{thm:2.6}, we conclude that $D$ is commutative and this is a
contradiction. 
\end{proof}

The following result follows immediately from Corollary \ref{cor:2.9}.

\begin{corollary}\label{cor:2.10}

If $D$ is a non-commutative division ring of type $2$, then $D^*/D'$ is not
finitely generated. 
\end{corollary}

\section{The radicality of subgroups}

In this section we study subgroups of $D^*$ which are radical over some
subring of $D$. To prove the next theorem we need the following useful
property of division rings of type $2$. 

\begin{lemma}\label{lem:3.1} Let $D$ be a division ring of type $2$ with
center $F$ and let $N$ be a subnormal subgroup of $D^*$. If for every pair
of elements $x, y\in N$, there exists some positive integer $n_{xy}$ such
that $x^{n_{xy}}y=yx^{n_{xy}}$, then $N\subseteq F$. 
\end{lemma}

\begin{proof}
Since $N$ is subnormal in $D^*$, there exists a series of subgroups
$$N=N_1\triangleleft N_2\triangleleft\ldots\triangleleft N_r=D^*.$$

Suppose that $x, y\in N$ and $K:=F(x, y)$. By putting 
$M_i=K\cap N_i,\, \forall i\in\{1, \ldots, r\},$ we obtain the following
series of subgroups:  
$$M_1\triangleleft M_2\triangleleft\ldots\triangleleft M_r=K^*.$$

For any $a\in M_1\leq N_1=N$, suppose that $n_{ax}$ and $n_{ay}$ are positive
integers such that 
$$a^{n_{ax}}x=xa^{n_{ax}} \mbox{ and } a^{n_{ay}}y=ya^{n_{ay}}.$$

Then, for $n:=n_{ax}n_{ay}$ we have
$$a^n=(a^{n_{ax}})^{n_{ay}}=(xa^{n_{ax}}x^{-1})^{n_{ay}}=xa^{n_{ax}n_{ay}}x^{-1}=xa^nx^{-1}$$
and
$$a^n=(a^{n_{ay}})^{n_{ax}}=(ya^{n_{ay}}y^{-1})^{n_{ax}}=ya^{n_{ay}n_{ay}}y^{-1}=ya^ny^{-1}.$$

Therefore $a^n\in Z(K)$. Hence $M_1$ is radical over $Z(K)$. By [6, Theorem
1], $M_1\subseteq Z(K)$. In particular, $x$ and $y$ commute with each
other. Consequently, $N$ is an abelian group. By [15, 14.4.4, p. 440],
$N\subseteq F$. 
\end{proof}

\begin{theorem}\label{thm:3.2}
Let $D$ be a division ring of type $2$ with center $F$, $K$ be a proper
division subring of $D$, and suppose that $N$ is a normal subgroup of
$D^*$. If $N$ is radical over $K$, then $N\subseteq F$. 
\end{theorem}

\begin{proof}
Suppose that $N$ is not contained in the center $F$. If 
$N\setminus K=\emptyset$, then $N\subseteq K$. By [15, p. 433], either
$K\subseteq F$ or $K=D$. Since $K\neq D$ by the assumption, it follows that
$K\subseteq F$. Hence $N\subseteq F$, which contradicts the assumption. Thus,
we have $N\setminus K\neq\emptyset$. 

Now, to complete the proof of our theorem we shall show that the elements of
$N$ satisfy the requirements of Lemma \ref{lem:3.1}. To this end, suppose that
$a, b\in N$. We examine the following cases: 

\noindent
{\em Case 1:}  $a\in K$.

{\em Subcase 1.1:} $b\not\in K$.

We shall prove that there exists some positive integer $n$ such that
$a^nb=ba^n$. Suppose that $a^nb\neq ba^n, \forall n\in\mathbb{N}$. Then,
$a+b\neq 0, a\neq \pm{1}$ and $b\neq \pm{1}$. So we have 
$$x=(a+b)a(a+b)^{-1}, y=(b+1)a(b+1)^{-1}\in N.$$

Since $N$ is radical over $K$, we can find positive integers $m_x$ and $m_y$
such that 
$$x^{m_x}=(a+b)a^{m_x}(a+b)^{-1}, y^{m_y}=(b+1)a^{m_y}(b+1)^{-1}\in K.$$

Putting $m=m_xm_y$, we have

$$x^m=(a+b)a^m(a+b)^{-1}, y^m=(b+1)a^m(b+1)^{-1}\in K.$$

Direct calculations give the equalities
$$x^mb-y^mb+x^ma-y^m=x^m(a+b)-y^m(b+1)=(a+b)a^m-(b+1)a^m=a^m(a-1),$$
from which we get
$$(x^m-y^m)b=a^m(a-1)+y^m-x^ma.$$

If $(x^m-y^m)\neq 0$, then $b=(x^m-y^m)^{-1}[a(a^m-1)+y^m-x^ma]\in K$,
contrary to the choice of $b$. Therefore $(x^m-y^m)= 0$ and
consequently, $a^m(a-1)=y^m(a-1)$. Since $a\neq 1,a^m=y^m=(b+1)a^m(b+1)^{-1}$
and it follows that $a^mb=ba^m$, a contradiction. 

{\em Subcase 1.2:} $b\in K$.

Consider an element $x\in N\setminus K$. Since $xb\not\in K$, by Subcase 1.1,
there exist positive integers $r, s$ such that 
$$a^rxb=xba^r \mbox{ and } a^sx=xa^s.$$

From these equalities it follows that
$$a^{rs}=(xb)^{-1}a^{rs}(xb)=b^{-1}(x^{-1}a^{rs}x)b=b^{-1}a^{rs}b,$$

and consequently, $a^{rs}b=ba^{rs}.$

\noindent
{\em Case 2:}  $a\not\in K$.

Since $N$ is radical over $K$, there exists some positive integer $m$ such
that $a^m\in K$.  By Case 1, there exists a positive integer $n$ such that
$a^{mn}b=ba^{mn}$. 
\end{proof}

Theorem \ref{thm:3.2} is closely related to the following conjecture of
Herstein in \cite{her2}: {\em ``For a division ring $D$, given a subnormal
subgroup $N$ of $D^*$. If $N$ is radical over the center $F$ of $D$, then
$N$ is central, i. e. $N\subseteq F$."} In Theorem \ref{thm:3.2}, the subgroup
$N$ is required to be radical over an arbitrary proper division subring $K$
of $D$, which does not necessarily coincide with the center $F$. Notice that
$N$ is required to be normal in $D^*$. So, the following question seems to be
interesting: {\em ``For a division ring $D$, given a subnormal subgroup $N$ of
$D^*$. Is $N$ contained in the center $F$of $D$, provided it is radical over
some proper division subring of $D$?"}   

Finally, we consider the question of the existence of maximal subgroups in $D$
which are radical over $F$. Recall that if $D$ is centrally finite of index
different from the characteristic of $F$, then $D^*$ contains no such
subgroups (see [1, Theorem 5]). Here, we consider the case when $D$ is of
type $2$ with $[D: F]=\infty$ and we prove that, if $char\,F=p > 0$, then $D^*$
contains no such subgroups. 

\begin{theorem}\label{thm:3.3}
Let $D$ be a division ring of type $2$ with center $F$ such that
$[D:F]=\infty$ and $char\, F = p>0$. Then the group $D^*$ contains no maximal
subgroups which are radical over $F$. 
\end{theorem}

\begin{proof}
Suppose that $M$ is a maximal subgroup of $D^*$ which is  radical over $F$. 
Put $G=D'\cap M$. For each $x\in G$, there exists a positive integer $n(x)$
such that $x^{n(x)}\in F$. It follows that $x^{n(x)}\in D'\cap F=Z(D')$. By
Theorem \ref{thm:2.3}, $Z(D')$ is torsion, so $x$ is periodic. Thus, $G$ is a
torsion group.  Since  $M'\leq G, M'$ is also torsion. For any $x,y\in M'$,
put $H=\ts{x,y}$ and $D_1=F(x,y)$. Then $n:=[D_1:F]<\infty$ and $H$ is a
torsion subgroup of $D_1^*\leq GL_n(F)$. By [12, (9.9'), p.  154], $H$ is
finite. Since $char F=p>0$, by [12, (13.3), p. 215], $H$ is cyclic. In
particular, $x$ and $y$ commute with each other, and consequently, $M'$ is
abelian. It follows that $M$ is a solvable group. Thus $M$ is a solvable
maximal  subgroup of $D^*$. By [1, Corollary 2, p. 432] and  [3, Theorem 6],
$[D:F]<\infty$, a contradiction. 
\end{proof}


\begin{thebibliography}{99}

\bibitem{Ak1} Akbari, S.;  Mahdavi-Hezavehi, M.; Mahmudi, M.G., Maximal
subgroups of $GL_1(D)$,  {\em  J. Algebra} 217 (1999), 422--433. 

\bibitem{Ak2} Akbari, S.;  Mahdavi-Hezavehi, M., Normal subgroups of $GL_n(D)$
are not finitely generated, {\em Proc. Amer. Math. Soc.} 128 (2000),
1627--1632. 

\bibitem{Ak3} Akbari, S.;  Ebrahimian, R.; Momenaee Kermani, H.; Salehi
Golsefidy, A.,  Maximal subgroups of $GL_n(D)$, {\em J. Algebra} 259 (2003),
201--225.  

\bibitem{dra} P. Draxl, {\em Skew Fields}, London Math. Soc. Lecture Note
Series 81 (1983), Cambridge Univ. Press. 

\bibitem{dren} V. Drensky, {\em Free Algebras and $PI$-Algebras}, Graduate
Course in Algebra, Springer-Verlag, Singapore, 2000. 

\bibitem{hai-huynh}  Bui Xuan Hai and Le Khac Huynh, On subgroups of the
multiplicative group of a division ring, {\em Vietnam J. Math.} 32 (2004),
21--24.  

\bibitem{hai-thin1} Bui Xuan Hai and Nguyen Van Thin,  On subnormal  and
maximal subgroups in division ring, {\em Southeast Asian Bull. Math.} 32
(2008), 931--937. 

\bibitem{hai-thin2} Bui Xuan  Hai and Nguyen Van  Thin,  On  locally nilpotent
subgroups of $GL_1(D)$, {\em Commun. Algebra} 37 (2009), 712--718. 

\bibitem{her} I. N. Herstein, {\em Noncommutative Rings}, The Carus
Mathematical Monographs. 

\bibitem{her2} I. N. Herstein,  Multiplicative commutators in division rings,
{\em Israel J. Math.}, 31 (1978), 180--188. 

\bibitem{kha} V. K. Kharchenko, Simple, prime and semiprime rings, {\em
Handbook of Algebra}, Vol. 1, North-Holland, Amsterdam, 1996, 761--812. 

\bibitem{lam} T.Y. Lam, {\em A First Course in Noncommutative Rings}, GTM 131
(1991), Springer-Verlag. 

\bibitem{mah} Mahdavi-Hezavehi M, Extension of valuations on derived groups of
division rings, {\em Commun. Algebra}, 23 (1995), 913--926. 

\bibitem{row} L. H. Rowen, {\em Ring Theory}, Vol. 1, Academic Press, New
York, 1988. 

\bibitem{scott} W. R. Scott, {\em Group Theory}, Dover Publication, INC, 1988.

\end{thebibliography}
\end{document}